# NONCENTRAL CONVERGENCE OF MULTIPLE INTEGRALS

By Ivan Nourdin and Giovanni Peccati

*Université Paris VI and Université Paris Ouest*

Fix $\nu > 0$, denote by $G(\nu/2)$ a Gamma random variable with parameter $\nu/2$ and let $n \geq 2$ be a fixed even integer. Consider a sequence $\{F_k\}_{k \geq 1}$ of square integrable random variables belonging to the $n$th Wiener chaos of a given Gaussian process and with variance converging to $2\nu$. As $k \to \infty$, we prove that $F_k$ converges in distribution to $2G(\nu/2) - \nu$ if and only if $E(F_k^4) - 12E(F_k^3) \to 12\nu^2 - 48\nu$.

**1. Introduction and main results.** Let $\mathfrak{H}$ be a real separable Hilbert space and, for $n \geq 1$, let $\mathfrak{H}^{\otimes n}$ (resp. $\mathfrak{H}^{\odot n}$) be the $n$th tensor product (resp. $n$th symmetric tensor product) of $\mathfrak{H}$. In what follows, we write

$$(1.1) \qquad X = \{X(h) : h \in \mathfrak{H}\}$$

to indicate a centered isonormal Gaussian process on $\mathfrak{H}$. For every $n \geq 1$, we denote by $I_n$ the isometry between $\mathfrak{H}^{\odot n}$ (equipped with the modified norm $\sqrt{n!} \|\cdot\|_{\mathfrak{H}^{\otimes n}}$) and the $n$th Wiener chaos of $X$. Note that, if $\mathfrak{H}$ is a $\sigma$-finite measure space with no atoms, then each random variable $I_n(h)$, $h \in \mathfrak{H}^{\odot n}$, has the form of a multiple Wiener–Itô integral of order $n$. For $n, m \geq 1$, $f \in \mathfrak{H}^{\odot n}$, $g \in \mathfrak{H}^{\odot m}$ and $p = 0, \ldots, n \wedge m$, we denote by $f \otimes_p g \in \mathfrak{H}^{\otimes(n+m-2p)}$ and $f \widetilde{\otimes}_p g \in \mathfrak{H}^{\odot(n+m-2p)}$, respectively, the $p$th contraction and the $p$th *symmetrized* contraction of $f$ and $g$ (a formal discussion of the properties of the previous objects is deferred to Section 2).

It is customary to call "Central Limit Theorem" (CLT in the sequel) any result describing the weak convergence of a (normalized) sequence of nonlinear functionals of $X$ toward a Gaussian law. Classic references for CLTs of this type are the works by Breuer and Major [1], Major [8], Giraitis and Surgailis [5] and Chambers and Slud [2]; the reader is also referred to the survey by Surgailis [14] and the references therein. More recently, Nualart and Peccati [11] proved the following result [here, and for the rest of the paper, we shall denote by $\mathscr{N}(0,1)$ the law of a Gaussian random variable with zero mean and unit variance].

THEOREM 1.1. *Fix an integer $n \geq 2$ and a sequence $\{f_k\}_{k \geq 1} \subset \mathfrak{H}^{\odot n}$ such that*

$$\lim_{k \to \infty} n! \|f_k\|_{\mathfrak{H}^{\otimes n}}^2 = \lim_{k \to \infty} E[I_n(f_k)^2] = 1. \tag{1.2}$$

*Then the following three conditions are equivalent:*

(i) $\lim_{k \to \infty} E[I_n(f_k)^4] = 3$;
(ii) *for every $p = 1, \ldots, n-1$, $\lim_{k \to \infty} \|f_k \otimes_p f_k\|_{\mathfrak{H}^{\otimes 2(n-p)}} = 0$;*
(iii) *as $k \to \infty$, the sequence $\{I_n(f_k)\}_{k \geq 1}$ converges in distribution to $N \sim \mathscr{N}(0,1)$.*

Theorem 1.1 is proved in Nualart and Peccati [11] by means of a stochastic calculus result, known as the *Dambis, Dubins and Schwarz theorem* (see, e.g., Revuz and Yor [13], Chapter V). In particular, Theorem 1.1 implies that the convergence in distribution of a sequence of multiple stochastic integrals toward a Gaussian random variable is completely determined by the asymptotic behavior of their second and fourth moments. As such, Theorem 1.1 can be seen as a drastic simplification of the classic "method of moments and diagrams" (see, for instance, the previously quoted works by Breuer, Major, Giraitis, Surgailis, Chambers and Slud).

The recent paper by Nualart and Ortiz-Latorre [10] contains a crucial methodological breakthrough, showing that one can prove Theorem 1.1 (as well as its multidimensional extensions) by using exclusively results from Malliavin calculus, such as integration by parts formulae and the duality properties of Malliavin derivatives and Skorohod integral operators. In particular, Nualart and Ortiz-Latorre prove that, for every $n \geq 2$ and for every sequence $\{I_n(f_k)\}_{k \geq 1}$ satisfying (1.2), either one of conditions (i)–(iii) in Theorem 1.1 is equivalent to the following: as $k \to \infty$,

$$\|D[I_n(f_k)]\|_{\mathfrak{H}}^2 \longrightarrow n \quad \text{in } L^2(\Omega), \tag{1.3}$$

where $D$ is the usual Malliavin derivative operator (see Section 2).

The principal aim of this paper is to prove several noncentral extensions of Theorem 1.1. Our main result is the following, which can be seen as a further simplification of the method of moments and diagrams, as applied to the framework of a non-Gaussian limit law. It should be compared with other noncentral limit theorems for nonlinear functionals of Gaussian fields, such as the ones proved by Taqqu [16, 17], Dobrushin and Major [3], Fox and Taqqu [4] and Terrin and Taqqu [18]; see also the survey by Surgailis [15] for further references in this direction.

THEOREM 1.2. *Let the previous notation prevail, fix $\nu > 0$ and let $F(\nu)$ be a real-valued random variable such that*

$$(1.4) \quad E(e^{i\lambda F(\nu)}) = \left(\frac{e^{-i\lambda}}{\sqrt{1-2i\lambda}}\right)^\nu, \qquad \lambda \in \mathbb{R}.$$

*Fix an* even *integer $n \geq 2$, and define*

$$(1.5) \quad c_n := \frac{1}{(n/2)!\binom{n-1}{n/2-1}^2} = \frac{4}{(n/2)!\binom{n}{n/2}^2}.$$

*Then for any sequence $\{f_k\}_{k\geq 1} \subset \mathfrak{H}^{\odot n}$ verifying*

$$(1.6) \quad \lim_{k\to\infty} n!\|f_k\|^2_{\mathfrak{H}^{\otimes n}} = \lim_{k\to\infty} E[I_n(f_k)^2] = 2\nu,$$

*the following six conditions are equivalent:*

(i) $\lim_{k\to\infty} E[I_n(f_k)^3] = E[F(\nu)^3] = 8\nu$ *and* $\lim_{k\to\infty} E[I_n(f_k)^4] = E[F(\nu)^4] = 12\nu^2 + 48\nu$;

(ii) $\lim_{k\to\infty} E[I_n(f_k)^4] - 12E[I_n(f_k)^3] = 12\nu^2 - 48\nu$;

(iii) $\lim_{k\to\infty} \|f_k \widetilde{\otimes}_{n/2} f_k - c_n \times f_k\|_{\mathfrak{H}^{\otimes n}} = 0$ *and* $\lim_{k\to\infty} \|f_k \widetilde{\otimes}_p f_k\|_{\mathfrak{H}^{\otimes 2(n-p)}} = 0$ *for every $p = 1, \ldots, n-1$ such that $p \neq n/2$;*

(iv) $\lim_{k\to\infty} \|f_k \widetilde{\otimes}_{n/2} f_k - c_n \times f_k\|_{\mathfrak{H}^{\otimes n}} = 0$ *and* $\lim_{k\to\infty} \|f_k \otimes_p f_k\|_{\mathfrak{H}^{\otimes 2(n-p)}} = 0$ *for every $p = 1, \ldots, n-1$ such that $p \neq n/2$;*

(v) *as $k \to \infty$, $\|D[I_n(f_k)]\|^2_{\mathfrak{H}} - 2nI_n(f_k) \longrightarrow 2n\nu$ in $L^2(\Omega)$, where $D$ is the Malliavin derivative operator;*

(vi) *as $k \to \infty$, the sequence $\{I_n(f_k)\}_{k\geq 1}$ converges in distribution to $F(\nu)$.*

REMARK 1.3.

1. The limit random variable $F(\nu)$ appearing in formula (1.4) is such that $F(\nu) \stackrel{\text{Law}}{=} 2G(\nu/2) - \nu$, where $G(\nu/2)$ has a Gamma law with parameter $\nu/2$, that is, $G(\nu/2)$ is a (a.s. strictly positive) random variable with density

$$g(x) = \frac{x^{\nu/2-1}e^{-x}}{\Gamma(\nu/2)}\mathbf{1}_{(0,\infty)}(x),$$

where $\Gamma$ is the usual Gamma function. Note that the following elementary relations have been implicitly used:

$$(1.7) \quad \begin{aligned} E(F(\nu)) &= 0, & E(F(\nu)^2) &= 2\nu, & E(F(\nu)^3) &= 8\nu, \\ E(F(\nu)^4) &= 12\nu^2 + 48\nu. \end{aligned}$$

2. When $\nu \geq 1$ is an integer, then $F(\nu)$ has a centered $\chi^2$ law with $\nu$ degrees of freedom. That is,

$$(1.8) \qquad F(\nu) \stackrel{\text{Law}}{=} \sum_{i=1}^{\nu} (N_i^2 - 1),$$

where $(N_1, \ldots, N_\nu)$ is a $\nu$-dimensional vector of i.i.d. $\mathcal{N}(0,1)$ random variables.

3. When $n \geq 1$ is an *odd* integer, there does not exist any sequence $\{I_n(f_k)\}_{k \geq 1}$, with $\{f_k\}_{k \geq 1} \subset \mathfrak{H}^{\odot n}$, such that $I_n(f_k)$ has bounded variances and $I_n(f_k)$ converges in distribution to $F(\nu)$ as $k \to \infty$. This is a consequence of the fact that any multiple integral of odd order has a third moment equal to zero, whereas $E(F(\nu)^3) = 8\nu > 0$.
4. The only difference between point (iii) and point (iv) of Theorem 1.2 is the symmetrization of the contractions of order $p \neq n/2$. One cannot dispense with the symmetrization of the contraction of order $n/2$. Note also that (iii) and (iv) do not depend on $\nu$; this means that, when applying either one of conditions (iii) and (iv), the dependence on $\nu$ is completely encoded by the normalization assumption (1.6).
5. In Proposition 4.1, we will use Theorem 1.1 in order to provide simple examples of sequences $\{I_n(f_k)\}_{k \geq 1}$ verifying both (1.6) and either one of the equivalent conditions (i)–(vi) of Theorem 1.2, for a given even integer $n \geq 4$ and a given *integer* $\nu \geq 1$.

Before going into details, we shall provide a short outline of the techniques used in the proof of Theorem 1.2. We will prove the following implications:

$$(\text{vi}) \to (\text{i}) \to (\text{ii}) \to (\text{iii}) \to (\text{iv}) \to (\text{v}) \to (\text{vi}).$$

The double implication (vi) → (i) → (ii) is trivial. The implication (ii) → (iii) is obtained by combining a standard version of the multiplication formula between multiple integrals with a result based on the integration by parts formulae of Malliavin calculus (see Lemma 2.1 below). The proof of (iii) → (iv) is purely combinatorial, whereas that of (iv) → (v) relies once again on multiplication formulae. Finally, to show (v) → (vi) we will adopt an approach similar to the one by Nualart and Ortiz-Latorre [10]. Our argument is as follows. Let us first observe that a sequence of random variables $\{I_n(f_k)\}_{k \geq 1}$ verifying (1.6) is tight and, therefore, by Prokhorov's theorem, it is relatively compact. As a consequence, to show the implication (v) → (vi), it is sufficient to prove that any subsequence $\{I_n(f_{k'})\}$, converging in distribution to some random variable $F_\infty$, must be necessarily such that $F_\infty \stackrel{\text{Law}}{=} F(\nu)$. This last property will be established by means of Malliavin calculus, by proving that condition (v) implies that the characteristic function $\phi_\infty$ of $F_\infty$ always solves the linear differential equation

$$(1.9) \qquad (1 - 2i\lambda)\phi'_\infty(\lambda) + 2\lambda\nu\phi_\infty(\lambda) = 0, \qquad \lambda \in \mathbb{R}, \phi_\infty(0) = 1.$$

Since the unique solution of (1.9) is given by the application $\lambda \mapsto E\{e^{i\lambda F(\nu)}\}$, the desired conclusion will follow immediately.

The paper is organized as follows. In Section 2, we present some preliminary results about Malliavin calculus. Section 3 contains the proof of Theorem 1.2 while, in Section 4, we give further refinements of Theorem 1.2.

**2. Preliminaries.** The reader is referred to the monograph by Nualart [9] for any unexplained notion or result discussed in this section. Let $\mathfrak{H}$ be a real separable Hilbert space. As in formula (1.1), we denote by $X$ an isonormal Gaussian process over $\mathfrak{H}$. Recall that, by definition, $X$ is a collection of centered and jointly Gaussian random variables indexed by the elements of $\mathfrak{H}$, defined on some probability space $(\Omega, \mathscr{F}, P)$ and such that, for every $h, g \in \mathfrak{H}$,

$$(2.1) \qquad E[X(h)X(g)] = \langle h, g \rangle_{\mathfrak{H}}.$$

We will systematically assume that $\mathscr{F}$ is generated by $X$. It is well known (see, e.g., Nualart [9], Chapter 1) that any random variable $F$ belonging to $L^2(\Omega, \mathscr{F}, P)$ admits the following chaotic expansion:

$$(2.2) \qquad F = \sum_{n=0}^{\infty} I_n(f_n),$$

where $I_0(f_0) := E[F]$, the series converges in $L^2(\Omega)$ and the kernels $f_n \in \mathfrak{H}^{\odot n}$, $n \geq 1$, are uniquely determined by $F$. Observe that $I_1(h) = X(h)$, $h \in \mathfrak{H}$, and that a random variable of the type $I_n(f)$, $f \in \mathfrak{H}^{\odot n}$, has finite moments of all orders (see, e.g., Janson [7], Chapter VI). As already pointed out, in the particular case where $\mathfrak{H} = L^2(A, \mathscr{A}, \mu)$, where $(A, \mathscr{A})$ is a measurable space and $\mu$ is a $\sigma$-finite and nonatomic measure, one has that $\mathfrak{H}^{\odot n} = L_s^2(A^n, \mathscr{A}^{\otimes n}, \mu^n)$ is the space of symmetric and square integrable functions on $A^n$. Moreover, for every $f \in \mathfrak{H}^{\odot n}$, $I_n(f)$ coincides with the multiple Wiener–Itô integral (of order $n$) of $f$ with respect to $X$ (see again Nualart [9], Chapter 1). For every $n \geq 0$, we write $J_n$ to indicate the orthogonal projection operator on the $n$th Wiener chaos associated with $X$. In particular, if $F \in L^2(\Omega, \mathscr{F}, P)$ is as in (2.2), then $J_n F = I_n(f_n)$ for every $n \geq 0$.

Let $\{e_k, k \geq 1\}$ be a complete orthonormal system in $\mathfrak{H}$. Given $f \in \mathfrak{H}^{\odot n}$ and $g \in \mathfrak{H}^{\odot m}$, for every $p = 0, \ldots, n \wedge m$, the $p$th contraction of $f$ and $g$ is the element of $\mathfrak{H}^{\otimes (n+m-2p)}$ defined as

$$(2.3) \quad f \otimes_p g = \sum_{i_1, \ldots, i_p = 1}^{\infty} \langle f, e_{i_1} \otimes \cdots \otimes e_{i_p} \rangle_{\mathfrak{H}^{\otimes p}} \otimes \langle g, e_{i_1} \otimes \cdots \otimes e_{i_p} \rangle_{\mathfrak{H}^{\otimes p}}.$$

Note that, in the particular case where $\mathfrak{H} = L^2(A, \mathscr{A}, \mu)$ (with $\mu$ nonatomic), one has that

$$(f \otimes_p g)(t_1, \ldots, t_{n+m-2p})$$

$$= \int_{A^p} f(t_1, \ldots, t_{n-p}, s_1, \ldots, s_p)$$
$$\times g(t_{n-p+1}, \ldots, t_{m+n-2p}, s_1, \ldots, s_p) \, d\mu(s_1) \cdots d\mu(s_p).$$

Moreover, $f \otimes_0 g = f \otimes g$ equals the tensor product of $f$ and $g$ while, for $n = m$, $f \otimes_n g = \langle f, g \rangle_{\mathfrak{H}^{\otimes n}}$. Note that, in general (and except for trivial cases), the contraction $f \otimes_p g$ is not a symmetric element of $\mathfrak{H}^{\otimes(n+m-2p)}$. As indicated in the Introduction, the canonical symmetrization of $f \otimes_p g$ is written $f \widetilde{\otimes}_p g$.

Let $\mathscr{S}$ be the set of all smooth cylindrical random variables of the form

$$F = g(X(\phi_1), \ldots, X(\phi_q)),$$

where $q \geq 1$, $g : \mathbb{R}^q \to \mathbb{R}$ is a smooth function with compact support and $\phi_i \in \mathfrak{H}$. The Malliavin derivative of $F$ with respect to $X$ is the element of $L^2(\Omega, \mathfrak{H})$ defined as

$$DF = \sum_{i=1}^{q} \frac{\partial g}{\partial x_i}(X(\phi_1), \ldots, X(\phi_q))\phi_i.$$

In particular, $DX(h) = h$ for every $h \in \mathfrak{H}$. By iteration, one can define the $m$th derivative $D^m F$ (which is an element of $L^2(\Omega, \mathfrak{H}^{\odot m})$) for every $m \geq 2$.

As usual, for $m \geq 1$, $\mathbb{D}^{m,2}$ denotes the closure of $\mathscr{S}$ with respect to the norm $\|\cdot\|_{m,2}$, defined by the relation

$$\|F\|_{m,2}^2 = E[F^2] + \sum_{i=1}^{m} E[\|D^m F\|_{\mathfrak{H}^{\otimes i}}^2].$$

The Malliavin derivative $D$ verifies the following *chain rule*: if $\varphi : \mathbb{R}^q \to \mathbb{R}$ is continuously differentiable with a bounded derivative and if $\{F_i\}_{i=1,\ldots,q}$ is a vector of elements of $\mathbb{D}^{1,2}$, then $\varphi(F_1, \ldots, F_q) \in \mathbb{D}^{1,2}$ and

$$D\varphi(F_1, \ldots, F_q) = \sum_{i=1}^{q} \frac{\partial \varphi}{\partial x_i}(F_1, \ldots, F_q) DF_i.$$

We denote by $\delta$ the adjoint of the operator $D$, also called the *divergence operator*. A random element $u \in L^2(\Omega, \mathfrak{H})$ belongs to the domain of $\delta$, noted $\mathrm{Dom}\,\delta$, if and only if it verifies

$$|E\langle DF, u \rangle_{\mathfrak{H}}| \leq c_u \sqrt{E(F^2)} \qquad \text{for any } F \in \mathscr{S},$$

where $c_u$ is a constant depending uniquely on $u$. If $u \in \mathrm{Dom}\,\delta$, then the random variable $\delta(u)$ is defined by the duality relationship (called "integration by parts formula"):

(2.4) $$E(F\delta(u)) = E\langle DF, u \rangle_{\mathfrak{H}},$$

which holds for every $F \in \mathbb{D}^{1,2}$. We will moreover need the following property: for every $F \in \mathbb{D}^{1,2}$ and every $u \in \mathrm{Dom}\,\delta$ such that $Fu$ and $F\delta(u) + \langle DF, u\rangle_{\mathfrak{H}}$ are square integrable, one has that $Fu \in \mathrm{Dom}\,\delta$ and

$$\delta(Fu) = F\delta(u) - \langle DF, u\rangle_{\mathfrak{H}}. \tag{2.5}$$

The operator $L$ is defined through the projection operators $\{J_n\}_{n\geq 0}$ as $L = \sum_{n=0}^{\infty} -nJ_n$, and is called the *infinitesimal generator of the Ornstein–Uhlenbeck semigroup*. It verifies the following crucial property: a random variable $F$ is an element of $\mathrm{Dom}\,L(=\mathbb{D}^{2,2})$ if and only if $F \in \mathrm{Dom}\,\delta D$ (i.e., $F \in \mathbb{D}^{1,2}$ and $DF \in \mathrm{Dom}\,\delta$), and in this case,

$$\delta DF = -LF.$$

Note that a random variable $F$ as in (2.2) is in $\mathbb{D}^{1,2}$ if and only if

$$\sum_{n=1}^{\infty} nn! \|f_n\|^2_{\mathfrak{H}^{\otimes n}} < \infty$$

and, in this case, $E[\|DF\|^2_{\mathfrak{H}}] = \sum_{n\geq 1} nn! \|f_n\|^2_{\mathfrak{H}^{\otimes n}}$. If $\mathfrak{H} = L^2(A, \mathscr{A}, \mu)$ (with $\mu$ nonatomic), then the derivative of a random variable $F$ as in (2.2) can be identified with the element of $L^2(A \times \Omega)$ given by

$$D_a F = \sum_{n=1}^{\infty} n I_{n-1}(f_n(\cdot, a)), \qquad a \in A. \tag{2.6}$$

The following lemma will be used in Section 3.

LEMMA 2.1. *Fix an integer $n \geq 2$ and set $F = I_n(f)$, with $f \in \mathfrak{H}^{\odot n}$. Then for every integer $s \geq 0$, we have*

$$E(F^s \|DF\|^2_{\mathfrak{H}}) = \frac{n}{s+1} E(F^{s+2}).$$

PROOF. We can write

$$E(F^s \|DF\|^2_{\mathfrak{H}}) = E(F^s \langle DF, DF\rangle_{\mathfrak{H}})$$

$$= \frac{1}{s+1} E(\langle DF, D(F^{s+1})\rangle_{\mathfrak{H}})$$

$$= \frac{1}{s+1} E(\delta DF \times F^{s+1}) \qquad \text{by integration by parts (2.4)}$$

$$= \frac{n}{s+1} E(F^{s+2})$$

by the property $\delta D = -L$ (which implies $\delta DF = nF$). □

**3. Proof of Theorem 1.2.** Throughout this section, $n \geq 2$ is an even integer, and $\{I_n(f_k)\}_{k \geq 1}$ is a sequence of multiple stochastic Wiener–Itô integrals of order $n$, such that condition (1.6) is satisfied for some $\nu > 0$.

3.1. *Proof of (vi) $\to$ (i) $\to$ (ii).* Since the sequence $\{I_n(f_k)\}_{k \geq 1}$ lives inside the $n$th chaos of $X$, and since condition (1.6) is in order, we deduce that, for every $p > 0$,

$$\sup_{k \geq 1} E[|I_n(f_k)|^p] < \infty \tag{3.1}$$

(see, e.g., Janson [7], Chapter V). This implies immediately that, if $\{I_n(f_k)\}_{k \geq 1}$ converges in distribution to $F(\nu)$, then, for every integer $p \geq 3$, $E(I_n(f_k)^p) \to E(F(\nu)^p)$. The implications (vi) $\to$ (i) $\to$ (ii) are therefore a direct consequence of (1.7).

3.2. *Proof of (ii) $\to$ (iii).* Suppose that (ii) holds. We start by observing that, due to the multiplication formulae between stochastic integrals (see Proposition 1.1.3 in Nualart [9]), we have

$$I_n(f_k)^2 = n! \|f_k\|_{\mathfrak{H}^{\otimes n}}^2 + \sum_{p=0}^{n-1} p! \binom{n}{p}^2 I_{2(n-p)}(f_k \widetilde{\otimes}_p f_k) \tag{3.2}$$

and

$$\|D[I_n(f_k)]\|_{\mathfrak{H}}^2 = nn! \|f_k\|_{\mathfrak{H}^{\otimes n}}^2 \tag{3.3}$$
$$+ n^2 \sum_{p=1}^{n-1} (p-1)! \binom{n-1}{p-1}^2 I_{2(n-p)}(f_k \widetilde{\otimes}_p f_k)$$

(see also Nualart and Ortiz-Latorre [10], Lemma 2). Relation (3.2) gives immediately that

$$E[I_n(f_k)^3] = n!(n/2)! \binom{n}{n/2}^2 \langle f_k, f_k \widetilde{\otimes}_{n/2} f_k \rangle_{\mathfrak{H}^{\otimes n}}. \tag{3.4}$$

On the other hand, we deduce from Lemma 2.1 (specialized to the case $s = 2$) that

$$E[I_n(f_k)^4] = \frac{3}{n} E[I_n(f_k)^2 \|D[I_n(f_k)]\|_{\mathfrak{H}}^2], \tag{3.5}$$

and, therefore, thanks to (3.2)–(3.3),

$$E[I_n(f_k)^4] = 3[n! \|f_k\|_{\mathfrak{H}^{\otimes n}}^2]^2 + \frac{3}{n} \sum_{p=1}^{n-1} n^2 (p-1)! \binom{n-1}{p-1}^2 p! \binom{n}{p}^2 \tag{3.6}$$
$$\times (2n-2p)! \|f_k \widetilde{\otimes}_p f_k\|_{\mathfrak{H}^{\otimes 2(n-p)}}^2.$$

In what follows, given two (deterministic) sequences $a(k)$ and $b(k)$, we write $a(k) \approx b(k)$ whenever $a(k) - b(k) \to 0$ as $k \to \infty$. Since (ii) and (1.6) hold, we deduce from (3.4)–(3.6) and condition (1.6), that

$$E[I_n(f_k)^4] - 12E[I_n(f_k)^3]$$

$$\approx [12\nu^2 - 48\nu] + \frac{3}{n}\sum_{\substack{p=1,\ldots,n-1 \\ p\neq n/2}} n^2(p-1)!\binom{n-1}{p-1}^2 p!\binom{n}{p}^2$$

$$\times (2n-2p)!\|f_k \widetilde{\otimes}_p f_k\|^2_{\mathfrak{H}^{\otimes 2(n-p)}}$$

(3.7) $\qquad + 24n!\|f_k\|^2_{\mathfrak{H}^{\otimes n}} + 3n(n/2-1)!\binom{n-1}{n/2-1}^2$

$$\times (n/2)!\binom{n}{n/2}^2 n!\|f_k \widetilde{\otimes}_{n/2} f_k\|^2_{\mathfrak{H}^{\otimes n}}$$

$$- 12n!(n/2)!\binom{n}{n/2}^2 \langle f_k, f_k \widetilde{\otimes}_{n/2} f_k\rangle_{\mathfrak{H}^{\otimes n}}.$$

Elementary simplifications give

$$24n!\|f_k\|^2_{\mathfrak{H}^{\otimes n}} + 3n(n/2-1)!\binom{n-1}{n/2-1}^2 (n/2)!\binom{n}{n/2}^2 n!\|f_k \widetilde{\otimes}_{n/2} f_k\|^2_{\mathfrak{H}^{\otimes n}}$$

$$- 12n!(n/2)!\binom{n}{n/2}^2 \langle f_k, f_k \widetilde{\otimes}_{n/2} f_k\rangle_{\mathfrak{H}^{\otimes n}}$$

$$= 24n!\|f_k\|^2_{\mathfrak{H}^{\otimes n}} + \frac{3}{2}(n!)^2 \binom{n}{n/2}^3 \|f_k \widetilde{\otimes}_{n/2} f_k\|^2_{\mathfrak{H}^{\otimes n}}$$

$$- 12n!(n/2)!\binom{n}{n/2}^2 \langle f_k, f_k \widetilde{\otimes}_{n/2} f_k\rangle_{\mathfrak{H}^{\otimes n}}$$

$$= \left\|2\sqrt{n!}\sqrt{6}f_k - \sqrt{\frac{3}{2}}\frac{(n!)^2\sqrt{n!}}{[(n/2)!]^3}f_k \widetilde{\otimes}_{n/2} f_k\right\|^2_{\mathfrak{H}^{\otimes n}}$$

$$= \frac{3}{2}\frac{(n!)^5}{[(n/2)!]^6}\|f_k \widetilde{\otimes}_{n/2} f_k - f_k \times c_n\|^2_{\mathfrak{H}^{\otimes n}},$$

where $c_n$ is defined in (1.5). This yields the desired conclusion.

3.3. *Proof of (iii) → (iv).* We can assume that $n \geq 4$. We shall introduce some further notation. Fix an integer $M \geq 1$, and denote by $\mathfrak{S}_{2M}$ the group of the $(2M)!$ permutations of the set $\{1,\ldots,2M\}$. We write $\pi_0$ to indicate the identity (trivial) permutation. Given a set $A$ and a vector $a = (a_1,\ldots,a_{2M}) \in A^{2M}$, for every $\pi \in \mathfrak{S}_{2M}$ we denote by $a_\pi = (a_{\pi(1)},\ldots,a_{\pi(2M)})$ the canonical action of $\pi$ on $a$. Note that, with this notation, one has $a = a_{\pi_0}$. For every $r = 0,\ldots,M$ and for $\pi,\sigma \in \mathfrak{S}_{2M}$, we write

$$\pi \sim_r \sigma$$

whenever the set $\{\pi(1),\ldots,\pi(M)\} \cap \{\sigma(1),\ldots,\sigma(M)\}$ contains exactly $r$ elements. For every $\pi \in \mathfrak{S}_{2M}$, there are exactly $M!^2 \binom{M}{r}^2$ permutations $\sigma$ such that $\pi \sim_r \sigma$. The implication (iii) $\to$ (iv) in the statement of Theorem 1.2 is a consequence of the following result.

PROPOSITION 3.1. *Let $n \geq 4$ be an even integer, and let $\{f_k\} \subset \mathfrak{H}^{\odot n}$ be a sequence of symmetric kernels. Then the following two conditions are equivalent:*

(A) $\|f_k \widetilde{\otimes}_p f_k\|_{\mathfrak{H}^{\otimes 2(n-p)}} \to 0$, $p = 1,\ldots,n-1$, $p \neq n/2$;
(B) $\|f_k \otimes_p f_k\|_{\mathfrak{H}^{\otimes 2(n-p)}} \to 0$, $p = 1,\ldots,n-1$, $p \neq n/2$.

PROOF. Since $\|f_k \otimes_p f_k\|_{\mathfrak{H}^{\otimes 2(n-p)}} \geq \|f_k \widetilde{\otimes}_p f_k\|_{\mathfrak{H}^{\otimes 2(n-p)}}$, the implication (B) $\Rightarrow$ (A) is trivial. Moreover, since

$$\|f_k \otimes_p f_k\|_{\mathfrak{H}^{\otimes 2(n-p)}} = \|f_k \otimes_{n-p} f_k\|_{\mathfrak{H}^{\otimes 2p}} \tag{3.8}$$

for every $p = 1,\ldots,n-1$, to show that (A) $\Rightarrow$ (B) it is sufficient to prove that (A) implies that $\forall p = 1,\ldots,\frac{n}{2}-1$,

$$\|f_k \otimes_p f_k\|_{\mathfrak{H}^{\otimes 2(n-p)}} \to 0. \tag{3.9}$$

Thanks to (3.8), and since $f_k \otimes_{n-1} f_k = f_k \widetilde{\otimes}_{n-1} f_k$, we immediately deduce that (A) implies that (3.9) holds for $p = 1$. This proves the implication (A) $\Rightarrow$ (B) in the case $n = 4$, so that from now on we can suppose that $n \geq 6$. The rest of the proof is done by recurrence. In particular, we shall show that, for every $q = 2,\ldots,\frac{n}{2}-1$, the following implication holds: if (A) is true and if (3.9) holds for $p = 1,\ldots,q-1$, then

$$\|f_k \otimes_q f_k\|_{\mathfrak{H}^{\otimes 2(n-q)}} \to 0.$$

Now fix $q = 2,\ldots,\frac{n}{2}-1$, suppose (A) is verified, and assume that (3.9) takes place for $p = 1,\ldots,q-1$. To simplify the discussion, we shall suppose (without loss of generality) that $\mathfrak{H} = L^2(A,\mathscr{A},\mu)$, where $\mu$ is $\sigma$-finite and nonatomic. Start by writing

$$\|f_k \widetilde{\otimes}_{n-q} f_k\|^2_{\mathfrak{H}^{\otimes 2q}} = \langle f_k \otimes_{n-q} f_k, f_k \widetilde{\otimes}_{n-q} f_k \rangle_{\mathfrak{H}^{\otimes 2q}}$$

$$= \frac{1}{(2q)!} \sum_{\pi \in \mathfrak{S}_{2q}} \int_{A^{2q}} f_k \otimes_{n-q} f_k(a_{\pi_0})$$

$$\times f_k \otimes_{n-q} f_k(a_\pi) \mu^{2q}(da).$$

Now, if $\pi \sim_0 \pi_0$ or $\pi \sim_q \pi_0$, one has that

$$\int_{A^{2q}} f_k \otimes_{n-q} f_k(a_{\pi_0}) \times f_k \otimes_{n-q} f_k(a_\pi) \mu^{2q}(da) = \|f_k \otimes_{n-q} f_k\|_{\mathfrak{H}^{\otimes 2q}}^2.$$

On the other hand, if $\pi \sim_p \pi_0$ for some $p = 1, \ldots, q-1$, then

$$(3.10) \quad \int_{A^{2q}} f_k \otimes_{n-q} f_k(a_{\pi_0}) \times f_k \otimes_{n-q} f_k(a_\pi) \mu^{2q}(da)$$

$$(3.11) \quad = \int_{A^{2(n-p)}} f_k \otimes_p f_k(a_{\pi^{[2(n-p)]}}) f_k \otimes_p f_k(a_{\sigma^{[2(n-p)]}}) \mu^{2(n-p)}(da),$$

where $(\pi^{[2(n-p)]}, \sigma^{[2(n-p)]}) \subset \mathfrak{S}_{2(n-p)}$ is any pair of permutations of $\{1, \ldots, 2(n-p)\}$ such that

$$\pi^{[2(n-p)]} \sim_{(q-p)} \sigma^{[2(n-p)]}.$$

Now, thanks to the recurrence assumption, and by Cauchy–Schwarz, we deduce that the expression in (3.11) converges to zero as $k \to \infty$, thus yielding that

$$0 = \lim_{k \to \infty} \|f_k \widetilde{\otimes}_{n-q} f_k\|_{\mathfrak{H}^{\otimes 2q}}^2 = \lim_{k \to \infty} \frac{2q!^2}{(2q)!} \|f_k \otimes_{n-q} f_k\|_{\mathfrak{H}^{\otimes 2q}}^2$$

$$= \lim_{k \to \infty} \frac{2}{\binom{2q}{q}} \|f_k \otimes_q f_k\|_{\mathfrak{H}^{\otimes 2(n-q)}}^2.$$

This concludes the proof. $\square$

3.4. *Proof of (iv) $\to$ (v)*. Suppose that (iv) holds. By using (3.3), we infer that $E[\|D[I_n(f_k)]\|_{\mathfrak{H}}^2] = nn! \|f_k\|_{\mathfrak{H}^{\otimes n}}^2 \to 2n\nu$. Moreover, by taking into account the orthogonality between multiple stochastic integrals of different orders and by using the multiplication formulae for multiple Wiener–Itô integrals, we have

$$E[I_n(f_k)\|D[I_n(f_k)]\|_{\mathfrak{H}}^2] = n^2(n/2-1)! \binom{n-1}{n/2-1}^2 n! \langle f_k \widetilde{\otimes}_{n/2} f_k, f_k \rangle_{\mathfrak{H}^{\otimes n}}$$

and

$$E[\|D[I_n(f_k)]\|_{\mathfrak{H}}^4] = n^4 \sum_{p=1}^n (p-1)!^2 \binom{n-1}{p-1}^4 (2n-2p)! \|f_k \widetilde{\otimes}_p f_k\|_{\mathfrak{H}^{\otimes 2(n-p)}}^2.$$

Now, define $c_n$ according to (1.5), and observe that (iv) and (1.6) imply that

$$\lim_{k \to \infty} \|f_k \widetilde{\otimes}_{n/2} f_k\|_{\mathfrak{H}^{\otimes n}}^2 = \lim_{k \to \infty} \|c_n f_k\|_{\mathfrak{H}^{\otimes n}}^2 = (2\nu c_n^2)/n!$$

and

$$\lim_{k \to \infty} \langle f_k \widetilde{\otimes}_{n/2} f_k, f_k \rangle_{\mathfrak{H}^{\otimes n}} = (2\nu c_n)/n!.$$

Thus, under (iv) one has that, as $k \to \infty$,

$$E(\|D[I_n(f_k)]\|^2 - 2nI_n(f_k) - 2n\nu)^2$$
$$= E[\|D[I_n(f_k)]\|^4] - 4nE[I_n(f_k)\|D[I_n(f_k)]\|^2]$$
$$+ 4n^2 E[I_n(f_k)^2] + 4n^2\nu^2 - 4n\nu E[\|D[I_n(f_k)]\|^2]$$
$$\longrightarrow 4\nu^2 n^2 + 2c_n^2 \nu n^4 (n/2-1)!^2 \binom{n-1}{n/2-1}^4$$
$$- 8c_n \nu n^3 (n/2-1)! \binom{n-1}{n/2-1}^2 + 8n^2\nu + 4n^2\nu^2 - 8n^2\nu^2 = 0.$$

3.5. *Proof of (v) → (vi).* Now we assume that (v) holds. We start by observing that condition (1.6) implies that the sequence of the laws of the random variables $\{I_n(f_k)\}_{k\geq 1}$ is tight (since it is bounded in $L^2(\Omega)$). By Prokhorov's theorem, we deduce that $\{I_n(f_k)\}_{k\geq 1}$ is relatively compact so that, to prove our claim, it is sufficient to show that any subsequence $\{I_n(f_{k'})\}$ converging in distribution to some random variable $F_\infty$ is necessarily such that $F_\infty \stackrel{\text{Law}}{=} F(\nu)$, where the law of $F(\nu)$ is defined by formula (1.4). From now on, and only for notational convenience, we assume that $\{I_n(f_k)\}$ itself converges to $F_\infty$. Also, for any $k \geq 1$, we let $\phi_k(\lambda) = E(e^{i\lambda I_n(f_k)})$ denote the characteristic function of $I_n(f_k)$, so that $\phi'_k(\lambda) = iE(I_n(f_k)e^{i\lambda I_n(f_k)})$. On the one hand, by the continuous mapping theorem, we have that

$$I_n(f_k)e^{i\lambda I_n(f_k)} \xrightarrow{\text{Law}} F_\infty e^{i\lambda F_\infty}.$$

Since boundedness in $L^2(\Omega)$ implies convergence of the expectations, we also deduce that $\phi'_k(\lambda) \to \phi'_\infty(\lambda)$ for any $\lambda \in \mathbb{R}$. On the other hand, we can write

$$\phi'_k(\lambda) = \frac{i}{n} E(\delta D[I_n(f_k)] \times e^{i\lambda I_n(f_k)}) \quad \text{since } \delta D = -L,$$
$$= \frac{i}{n} E(\langle D[I_n(f_k)], D(e^{i\lambda I_n(f_k)})\rangle_{\mathfrak{H}}) \quad \text{by integration by parts (2.4)},$$
$$= -\frac{\lambda}{n} E(e^{i\lambda I_n(f_k)} \|D[I_n(f_k)]\|_{\mathfrak{H}}^2).$$

Since (v) is in order, we deduce that, as $k \to \infty$,

$$\phi'_k(\lambda) + 2\lambda E(e^{i\lambda I_n(f_k)} I_n(f_k)) + 2\lambda \nu E(e^{i\lambda I_n(f_k)}) \to 0.$$

As a consequence, $\phi_\infty$ must necessarily solve the linear differential equation (1.9), yielding

$$\phi_\infty(\lambda) = \left(\frac{e^{-i\lambda}}{\sqrt{1-i2\lambda}}\right)^\nu = E(e^{i\lambda F(\nu)}), \quad \lambda \in \mathbb{R}.$$

This concludes the proof of Theorem 1.2.

**4. Further remarks.** When $\nu \geq 1$ is an integer, one can use Theorem 1.1 in order to obtain examples of sequences of multiple integrals $\{I_{2m}(f_k)\}_{k\geq 1}$ ($m \geq 2$ fixed) satisfying either one of conditions (i)–(vi) in Theorem 1.2. This fact is summarized in the following statement.

PROPOSITION 4.1. *Let $m \geq 2$ and $\nu \geq 1$ be integers and, for $i = 1, \ldots, \nu$, let $\{g_k^i\}_{k\geq 1} \subset \mathfrak{H}^{\odot m}$ be a sequence of kernels such that, as $k \to \infty$: (i) $m!\langle g_k^i, g_k^j \rangle_{\mathfrak{H}^{\otimes m}} \to \delta_{ij}$ for every $1 \leq i, j \leq \nu$ ($\delta_{ij}$ stands for the Kronecker symbol), (ii) $\|g_k^i \otimes_p g_k^i\|_{\mathfrak{H}^{\otimes 2(m-p)}} \to 0$ for every $1 \leq i \leq \nu$ and $1 \leq p \leq m - 1$. Then the sequence $\{I_{2m}(f_k)\}_{k\geq 1}$, where $f_k = \sum_{i=1}^{\nu} g_k^i \widetilde{\otimes} g_k^i \in \mathfrak{H}^{\odot 2m}$, converges in distribution to $F(\nu) \stackrel{\text{Law}}{=} \sum_{i=1}^{\nu}(N_i^2 - 1)$, where $(N_1, \ldots, N_\nu)$ is a vector of i.i.d. $\mathcal{N}(0,1)$ random variables.*

PROOF. Since conditions (i) and (ii) are in order, we deduce from [12] that
$$(I_m(g_k^1), \ldots, I_m(g_k^\nu)) \stackrel{\text{Law}}{\longrightarrow} \mathcal{N}_\nu(0, \text{Id}),$$
where $\mathcal{N}_\nu(0, \text{Id})$ stands for a $\nu$-dimensional Gaussian vector with zero mean and covariance equal to the identity matrix. On the other hand, as a consequence of the multiplication formula for Wiener–Itô integrals, we have
$$\sum_{i=1}^{\nu} I_m(g_k^i)^2 = \sum_{i=1}^{\nu} m!\|g_k^i\|_{\mathfrak{H}^{\otimes m}}^2 + \sum_{i=1}^{\nu}\sum_{p=1}^{m-1} p!\binom{m}{p}^2 I_{2(m-p)}(g_k^i \widetilde{\otimes}_p g_k^i)$$
$$+ \sum_{i=1}^{\nu} I_{2m}(g_k^i \widetilde{\otimes} g_k^i).$$

Since $\sum_{i=1}^{\nu} I_{2m}(g_k^i \widetilde{\otimes} g_k^i) = I_{2m}(f_k)$ (by linearity) and $\|g_k^i \widetilde{\otimes}_p g_k^i\|_{\mathfrak{H}^{\otimes 2(m-p)}} \leq \|g_k^i \otimes_p g_k^i\|_{\mathfrak{H}^{\otimes 2(m-p)}}$, the conclusion is immediately obtained. □

A refinement of Theorem 1.2 is the following.

PROPOSITION 4.2. *Let $n \geq 4$ be an even integer, and let $\{I_n(f_k)\}_{k\geq 1}$ be a sequence of multiple integrals verifying (1.6) and satisfying either one of conditions (i)–(vi) of Theorem 1.2. Then for every $h_1, \ldots, h_r \in \mathfrak{H}$ ($r \geq 1$), the vector $(I_n(f_k), I_1(h_1), \ldots, I_1(h_r))$ converges in law to $(F(\nu), I_1(h_1), \ldots, I_1(h_r))$ as $k \to \infty$, where $F(\nu) \stackrel{\text{Law}}{=} 2G(\nu/2) - \nu$ is independent of $X$.*

PROOF. By the definitions of the contractions of order 1 and $n-1$, for every fixed $j \in \{1, \ldots, r\}$, one has that
$$\|f_k \otimes_1 h_j\|_{\mathfrak{H}^{\otimes(n-1)}}^2 = \langle f_k \otimes_{n-1} f_k, h_j \otimes h_j\rangle_{\mathfrak{H}^{\otimes 2}} \leq \|f_k \otimes_{n-1} f_k\|_{\mathfrak{H}^{\otimes 2}} \|h_j\|_{\mathfrak{H}}^2 \underset{k\to\infty}{\longrightarrow} 0,$$

where the last convergence is a consequence of point (iv) in Theorem 1.2, and of the fact that $n \geq 4$. On the other hand,

$$E\langle D[I_n(f_k)], h_j\rangle_{\mathfrak{H}}^2 = nn! \|f_k \widetilde{\otimes}_1 h_j\|_{\mathfrak{H}^{\otimes(n-1)}}^2 \underset{k\to\infty}{\longrightarrow} 0$$

for any fixed $j$, so that one can finish the proof by simply mimicking the arguments displayed in Nualart and Ortiz-Latorre [10], proof of Theorem 7. □

REMARK 4.3.

1. When $n = 2$, the statement of Proposition 4.2 is not true in general. As a counterexample, one can consider a constant sequence $I_2(f_k)$, $k \geq 1$, such that $f_k = h \otimes h$ for every $k$, and $\|h\|_{\mathfrak{H}} = 1$.
2. Proposition 4.2 can be reformulated by saying that $I_n(f_k)$ *converges $\sigma\{X\}$-stably* to $F(\nu)$ (see, e.g., Jacod and Shiryayev [6] for an exhaustive discussion of stable convergence).

Proposition 4.2 yields a refinement of a well-known result (see, e.g., Janson [7], Chapter VI, Corollary 6.13), stating that Wiener chaoses of order $n > 2$ do not contain any Gamma random variable (our refinement consists in a further restriction on moments).

COROLLARY 4.4. *Fix a real $\nu > 0$ and an even integer $n \geq 4$. Let $I_n(f)$ be such that $E(I_n(f)^2) = 2\nu$. Then $I_n(f)$ cannot be equal in law to $2G(\nu/2) - \nu$, where $G(\nu/2)$ stands for a Gamma random variable of parameter $\nu/2$, and $E(I_n(f)^4) - 12E(I_n(f)^3) > 12\nu^2 - 48\nu$.*

PROOF. According to Proposition 4.2, if $I_n(f)$ was equal in law to $2G(\nu/2) - \nu$ (or if $E(I_n(f)^4) - 12E(I_n(f)^3) = 12\nu^2 - 48\nu$), then $I_n(f)$ would be independent of $X$. Plainly, this is only possible if $f = 0$, which is absurd, since $\|f\|_{\mathfrak{H}^{\otimes n}}^2 = 2\nu/n!$. The fact that $E(I_n(f)^4) - 12E(I_n(f)^3)$ cannot be less than $12\nu^2 - 48\nu$ derives from a straightforward modification of the calculations following formula (3.7). □

The following result characterizes the stable convergence of double integrals. The proof (omitted) is analogous to that of Proposition 4.2.

PROPOSITION 4.5. *Fix $\nu > 0$. Let the sequence $\{I_2(f_k)\}_{k\geq 1}$ be such that $E(I_2(f_k)^2) \to 2\nu$ and either one of conditions (i)–(vi) of Theorem 1.2 are satisfied. Then for every $h_1, \ldots, h_r \in \mathfrak{H}$ ($r \geq 1$), the vector $(I_2(f_k), I_1(h_1), \ldots, I_1(h_r))$ converges in law to $(F(\nu), I_1(h_1), \ldots, I_1(h_r))$, where $F(\nu) \stackrel{\text{Law}}{=} 2G(\nu/2) - \nu$ is independent of $X$, if and only if for any $j \in \{1, \ldots, r\}$,*

$$(4.1) \qquad \langle f_k \otimes_1 f_k, h_j \otimes h_j\rangle_{\mathfrak{H}^{\otimes 2}} \to 0.$$

*In particular, $I_2(f_k)$ is asymptotically independent of $X$ if and only if (4.1) is verified for any $j$.*

**Acknowledgments.** We are grateful to S. Kwapień, M. S. Taqqu and M. Yor for inspiring discussions on the subject of this paper. We also thank an anonymous referee for valuable suggestions.

## REFERENCES

[1] BREUER, P. and MAJOR, P. (1983). Central limit theorems for nonlinear functionals of Gaussian fields. *J. Multivariate Anal.* **13** 425–441. MR716933
[2] CHAMBERS, D. and SLUD, E. (1989). Central limit theorems for nonlinear functionals of stationary Gaussian processes. *Probab. Theory Related Fields* **80** 323–346. MR976529
[3] DOBRUSHIN, R. L. and MAJOR, P. (1979). Non-central limit theorems for nonlinear functionals of Gaussian fields. *Z. Wahrsch. Verw. Gebiete* **50** 27–52. MR550122
[4] FOX, R. and TAQQU, M. S. (1985). Noncentral limit theorems for quadratic forms in random variables having long-range dependence. *Ann. Probab.* **13** 428–446. MR781415
[5] GIRAITIS, L. and SURGAILIS, D. (1985). CLT and other limit theorems for functionals of Gaussian processes. *Z. Wahrsch. Verw. Gebiete* **70** 191–212. MR799146
[6] JACOD, J. and SHIRYAEV, A. N. (1987). *Limit Theorems for Stochastic Processes. Grundlehren der Mathematischen Wissenschaften [Fundamental Principles of Mathematical Sciences]* **288**. Springer, Berlin. MR959133
[7] JANSON, S. (1997). *Gaussian Hilbert Spaces. Cambridge Tracts in Mathematics* **129**. Cambridge Univ. Press, Cambridge. MR1474726
[8] MAJOR, P. (1981). *Multiple Wiener–Itô Integrals: With Applications to Limit Theorems. Lecture Notes in Math.* **849**. Springer, Berlin. MR611334
[9] NUALART, D. (2006). *The Malliavin Calculus and Related Topics*, 2nd ed. Springer, Berlin. MR2200233
[10] NUALART, D. and ORTIZ-LATORRE, S. (2008). Central limit theorems for multiple stochastic integrals and Malliavin calculus. *Stochastic Process. Appl.* **118** 614–628. MR2394845
[11] NUALART, D. and PECCATI, G. (2005). Central limit theorems for sequences of multiple stochastic integrals. *Ann. Probab.* **33** 177–193. MR2118863
[12] PECCATI, G. and TUDOR, C. A. (2005). Gaussian limits for vector-valued multiple stochastic integrals. In *Séminaire de Probabilités XXXVIII. Lecture Notes in Math.* **1857** 247–262. Springer, Berlin. MR2126978
[13] REVUZ, D. and YOR, M. (1999). *Continuous Martingales and Brownian Motion*, 3rd ed. *Grundlehren der Mathematischen Wissenschaften [Fundamental Principles of Mathematical Sciences]* **293**. Springer, Berlin. MR1725357
[14] SURGAILIS, D. (2003). CLTs for polynomials of linear sequences: Diagram formula with illustrations. In *Theory and Applications of Long-range Dependence* 111–127. Birkhäuser Boston, Boston, MA. MR1956046
[15] SURGAILIS, D. (2003). Non-CLTs: $U$-statistics, multinomial formula and approximations of multiple Itô–Wiener integrals. In *Theory and Applications of Long-range Dependence* 129–142. Birkhäuser, Boston. MR1956047
[16] TAQQU, M. S. (1974/75). Weak convergence to fractional Brownian motion and to the Rosenblatt process. *Z. Wahrsch. Verw. Gebiete* **31** 287–302. MR0400329

[17] Taqqu, M. S. (1979). Convergence of integrated processes of arbitrary Hermite rank. *Z. Wahrsch. Verw. Gebiete* **50** 53–83. MR550123
[18] Terrin, N. and Taqqu, M. S. (1990). A noncentral limit theorem for quadratic forms of Gaussian stationary sequences. *J. Theoret. Probab.* **3** 449–475. MR1057525

Laboratoire de Probabilités
 et Modèles Aléatoires
Université Paris VI
Boîte courrier 188, 4 Place Jussieu
75252 Paris Cedex 5
France
E-mail: ivan.nourdin@upmc.fr

Equipe Modal'X
Université Paris Ouest — Nanterre la Défense
frm-e00 Avenue de la République
92000 Nanterre
and
LSTA Université Paris VI
France
E-mail: giovanni.peccati@gmail.com